\documentclass{article}

\usepackage{amssymb}
\usepackage{amsmath}
\usepackage{tensor}

\newtheorem{theorem}{Theorem}[section]

\begin{document}

\parindent0pt
\title{Every projective variety is a quiver Grassmannian}
\author{Markus Reineke}
\maketitle

Quiver Grassmannians are varieties ${\rm Gr}_{\bf e}^Q(V)$ parametrizing subrepresentations of dimension vector ${\bf e}$ of a representation $V$ of a quiver $Q$ \cite{SchofieldGeneric}. The aim of this note is to prove the claim in the title, which implies that no special properties of quiver Grassmannians can be expected without restricting the choices of $Q$, $V$ and ${\bf e}$. More precisely we prove:

\begin{theorem} Every projective variety is isomorphic to a quiver Grassmannian ${\rm Gr}^Q_{\bf e}(V)$, for an acyclic quiver $Q$ with at most three vertices, a Schurian representation $V$, and a thin dimension vector ${\bf e}$ (i.e.~${\bf e}_i\leq 1$ for all $i\in Q_0$).
\end{theorem}

The author would like to thank B. Keller for posing the question answered here, and G. Cerulli-Irelli, E. Feigin, O. Lorscheid and A. Zelevinsky for helpful discussions.\\[2ex]
So let $X$ be an arbitrary projective variety, given by a set of homogeneous polynomial equations on the coordinates of ${\bf P}^n$. Without loss of generality, we can assume these equations to be homogeneous of the same degree $d$. Using the $d$-uple embedding $j:{\bf P}^n\rightarrow{\bf P}^{M-1}$, we can define $j(X)\simeq X$ by equations defining $j({\bf P}^n)$ inside ${\bf P}^{M-1}$, together with linear equations $\varphi_1,\ldots,\varphi_k$ on the coordinates of ${\bf P}^{M-1}$.\\

To make this more explicit, let $M_{n,d}$ be the set of tuples ${\bf m}=(m_0,\ldots,m_n)$ in ${\bf N}^{n+1}$ summing up to $d$, so that $M$ is the cardinality of $M_{n,d}$, and $j$ maps homogeneous coordinates $(x_0:\ldots:x_n)$ to $(:x^{\bf m}:)_{{\bf m}\in M_{n,d}}$, where $x^{\bf m}=x_0^{m_0}\ldots x_n^{m_n}$. A point $x$ in ${\bf P}^M$ with homogeneous coordinates $(:x_{\bf m}:)_{{\bf m}\in M_{n,d}}$ belongs to the image of $j$ if and only if the quadratic conditions $x_{\bf m}x_{\bf m'}=x_{\bf m''}x_{\bf m'''}$ for all ${\bf m}+{\bf m'}={\bf m''}+{\bf m'''}$ are fulfilled. Actually it suffices to impose the conditions (denoting by ${\bf e}_i$ the tuple with entry $1$ in the $i$-th position)
$$x_{{\bf n}+{\bf e}_i}x_{{\bf n'}+{\bf e}_j}=x_{{\bf n}+{\bf e}_j}x_{{\bf n'}+{\bf e}_i}\mbox{ for all }{\bf n}\in M_{n,d-1},\, i=0,\ldots,n.$$
These conditions can be rewritten as follows: we define a matrix $A(x)$ with rows indexed by the ${\bf n}\in M_{n,d-1}$ and columns indexed by the $i=0,\ldots,n$, with $({\bf n},i)$-th entry being $x_{{\bf n}+{\bf e}_i}$. Then the previous quadratic conditions can be compactly reformulated as the rank of $A$ being $1$, because they translate precisely to vanishing of all two-by-two minors of $A$. Thus we can define (a variety isomorphic to) $X$ as the set of points $x\in{\bf P}^M$ such that $A(x)$ has rank $1$ and $\varphi_j(x)=0$ for $j=1,\ldots,k$.\\

These defining equations are easily realized as a quiver Grassmannian: let $Q$ be the quiver with set of vertices $\{1,2,3\}$, with $k$ arrows from $2$ to $1$ and $n+1$ arrows from $2$ to $3$. Define a representation $V$ of $Q$ of dimension vector $(1,M,M')$, for $M'$ the cardinality of $M_{n,d-1}$, as follows: let $V_1={\bf C}$, the space $V_2$ has a basis $(v_{\bf m})_{{\bf m}\in M_{n,d}}$, and $V_3$ has a basis $(v_{\bf n})_{{\bf n}\in M_{n,d-1}}$. The maps representing the arrows from $2$ to $1$ are given by the linear forms $\varphi_j$. The $i$-th arrow from $2$ to $3$, with $i=0,\ldots,n$, is given by the linear map $f_i$ sending $v_{\bf m}$ to $v_{{\bf m}-{\bf e}_i}$ if $m_i\not=0$, and to $0$ otherwise. Define the dimension vector ${\bf e}=(0,1,1)$. Then a point of ${\rm Gr}_{\bf e}(V)$ is determined by a vector $v$ in $V_2$ which is annihilated by all $\varphi_j$, and such that the span of the images $f_0(v),\ldots,f_n(v)$ is one-dimensional. Writing $v=\sum_{{\bf m}\in M_{n,d}}x_{\bf m}v_{\bf m}$, the image $f_i(v)$, written in coordinates $v_{\bf n}$, is precisely the $i$-th column of $A(x)$, so that the span being one-dimensional means that the rank of $A(x)$ equals $1$. This proves all claims except $V$ being Schurian.\\

To prove this, assume we are given endomorphisms $\varphi$ of $V_2$ and $\psi$ of $V_3$ such that $f_i\varphi=\psi f_i$ for all $i=0,\ldots,n$. Writing $\varphi(v_{\bf m})=\sum_{{\bf m}'}a_{{\bf m},{\bf m}'}v_{{\bf m}'}$ and $\psi(v_{\bf n})=\sum_{{\bf n}'}b_{{\bf n},{\bf n}'}v_{{\bf n}'}$, this condition means that, for all $i=0,\ldots,n$ such that ${m}'_i\not=0$, we have
\begin{description}
\item[(i)] $a_{{\bf m},{\bf m}'}=0$ if ${m}_i=0$,
\item[(ii)] $a_{{\bf m},{\bf m}'}=b_{{\bf m}-{\bf e}_i,{\bf m}'-{\bf e}_i}$ if ${m}_i\not=0$.
\end{description}
The second condition implies that $a_{{\bf m},{\bf m}'}=a_{{\bf m}-{\bf e}_i+{\bf e}_j,{\bf m}'-{\bf e}_i+{\bf e}_j}$ if ${ m}_i\not=0\not={ m}_i$, and thus the following translation principle:  $a_{{\bf m},{\bf m}'}=a_{{\bf m}+{\bf q},{\bf m}'+{\bf q}}$ for all ${\bf q}=(q_0,\ldots,q_n)$ with $\sum_iq_i=0$ whenever all coefficients are defined. Given a pair ${\bf m}\not={\bf m}'$, we can find an index $i$ such that ${m}_i<{m}'_i$, thus $a_{{\bf m},{\bf m}'}=a_{{\bf m}-{m}_i{\bf e}_i+{m}_i{\bf e}_j,{\bf m}'-{ m}_i{\bf e}_i+{m}_i{\bf e}_j}$ for some $j\not=i$ by translation. Condition (i) now implies $a_{{\bf m},{\bf m}'}=0$. Moreover, all coefficients $a_{{\bf m},{\bf m}'}$ coincide by translation, thus $a_{{\bf m},{\bf m}'}=C\delta_{{\bf m},{\bf m}'}$ for some scalar $C$. All coefficients $b_{{\bf n},{\bf n}'}$ being determined by the $a_{{\bf m},{\bf m}'}$ by condition (ii), we have $b_{{\bf n},{\bf n}'}=C\delta_{{\bf n},{\bf n}'}$ as well. Now a non-zero $\varphi_k$ forces the endomorphism of $V$ to be just $C$ times the identity. In case all $\varphi_k$ are zero, we just omit the vertex $1$ from the quiver. The theorem is proved.




\begin{thebibliography}{9}
\bibitem{SchofieldGeneric} A.~Schofield, \emph{General representations of quivers}, Proc. London Math.
  Soc. (3) \textbf{65} (1992), no.~1, 46--64. 
\end{thebibliography}
\end{document}